\documentclass{amsart}

\usepackage{natbib}
\usepackage{latexcad}
\usepackage{amssymb}

\theoremstyle{definition}
\newtheorem{definition}{Definition}[section]
\newtheorem{example}[definition]{Example}
\newtheorem{remark}[definition]{Remark}

\theoremstyle{plain}
\newtheorem{lemma}[definition]{Lemma}
\newtheorem{theorem}[definition]{Theorem}
\newtheorem{proposition}[definition]{Proposition}
\newtheorem{corollary}[definition]{Corollary}
\newtheorem{conjecture}[definition]{Conjecture}

\def\idname#1#2#3#4{#1{}_{#2}{#4}{}_{#3}}
\def\sidname#1#2#3{{}_{#1}{#3}{}_{#2}}
\def\wreath#1#2{#1\ \mathrm{wr}_X\ #2}

\title[Linear Groupoids and Wreath Products]
{Linear Groupoids and the Associated Wreath Products}

\author{J.~D.~Phillips}

\address{Department of Mathematics \& Computer Science, Wabash College,
Crawfordsville, Indiana 47933, U.S.A.}

\email{phillipj@wabash.edu}

\author{Petr Vojt\v{e}chovsk\'y}

\address{Department of Mathematics, University of Denver, 2360 S Gaylord St,
Denver, CO, 80208, U.S.A.}

\email{petr@math.du.edu}

\begin{document}

\begin{abstract}
A groupoid identity is said to be linear of length $2k$ if the same $k$
variables appear on both sides of the identity exactly once. We classify and
count all varieties of groupoids defined by a single linear identity. For
$k=3$, there are $14$ nontrivial varieties and they are in the most general
position with respect to inclusion. Hentzel et. al. \cite{HJM1993} showed
that the linear identity $(xy)z = y(zx)$ implies commutativity and
associativity in all products of at least 5 factors. We complete their
project by showing that no other linear identity of any length behaves this
way, and by showing how the identity $(xy)z = y(zx)$ affects products of
fewer than 5 factors; we include distinguishing examples produced by the
finite model builder Mace4. The corresponding combinatorial results for
labelled binary trees are given. We associate a certain wreath product with
any linear identity. Questions about linear groupoids can therefore be
transferred to groups and attacked by group-theoretical computational tools,
e.g., GAP. Systematic notation and diagrams for linear identities are
devised. A short equational basis for Boolean algebras involving the identity
$(xy)z = y(zx)$ is presented, together with a proof produced by the automated
theorem prover Otter.
\end{abstract}

\keywords{linear groupoid, linear identity, balanced identity, strictly
balanced identity, the identity $(xy)z=y(zx)$, binary tree, wreath product,
Robbins axiom, boolean algebra, identity-hedron}

\subjclass{Primary: 20N0., Secondary: 18B40, 20B40, 20N05.}

\maketitle

\section{Motivation}\label{Sc:M}

\noindent It is customary to call an identity \emph{balanced} if the same
variables occur on both sides of the identity the same number of times. When
each of the $k$ variables of a balanced identity $\iota$ appears on each side
of $\iota$ exactly once, $\iota$ is called \emph{strictly balanced} or
\emph{linear of length} $2k$. We use the name \emph{linear} in this paper.

Thus, the \emph{associative law} $x(yz)=(xy)z$ is a linear identity of length
$6$, and the \emph{medial law} $(xy)(uv) = (xu)(yv)$ is a linear identity of
length $8$.

There does not seem to be any systematic account of groupoids satisfying a
linear identity, although several specific identities have been studied in
considerable detail. For instance, Je\v{z}ek and Kepka wrote a series of
papers on linear identities with identical bracketings on both sides, e.g.,
the \emph{medial groupoids} defined by the above medial law \cite{JKMed}, the
\emph{left} (resp. \emph{right}) \emph{permutable groupoids} defined by
$x(yz)=x(zy)$ (resp.\ $(xy)z=(xz)y$) \cite{JKPer}, and the \emph{left}
(resp.\ \emph{right}) \emph{modular groupoids} defined by $x(yz)=z(yx)$
(resp.\ $(xy)z=(zy)x$) \cite{JKMod}. These papers deal mostly with a
representation of linear groupoids by means of commutative semigroups, with
the description of all (finite) simple linear groupoids in a given variety,
and with universal algebraic properties of the varieties of linear groupoids.

We were drawn to the subject by the fascinating identity
\begin{equation}\label{Eq:H}
    (xy)z=y(zx),
\end{equation}
which, as far as we know, has not been named yet. Hentzel, Jacobs and Muddana
\cite{HJM1993} showed that for any groupoid $G$ satisfying \eqref{Eq:H} and
for any product of $m\ge 5$ elements of $G$, the $m$ factors commute and
associate, i.e., the result of the product is independent of parentheses and
of the order in which the elements are multiplied. This sounds paradoxical,
since it is certainly not true for $m=3$, and one would intuitively expect
the situation to become more complex with increasing $m$.

No explanation (beside a proof!) for this phenomenon is offered in
\cite{HJM1993}. A superficial explanation could go as follows: the longer the
products become, the more ways there are in which the substitution rule
\eqref{Eq:H} can be applied to them. Unfortunately, it is not clear at all why
this should overpower the growing number of possible products, or why it only
works for \eqref{Eq:H} and not for other linear identities.

\subsection{Contents}
We introduce a systematic notation for linear identities, and capture the
behavior of linear identities as substitutions in diagrams called
identity-hedrons. Given two linear identities, we decide when one implies the
other. Consequently, we can count how many distinct varieties of groupoids
defined by a single linear identity of given length there are. The answer
depends on the number of cyclic subgroups of symmetric groups. We show that
the only linear identity that implies associativity and commutativity in
sufficiently long products is $(\ref{Eq:H})$. This result can be restated in
terms of transformations of labelled binary trees. We introduce a canonical
way of constructing a certain subgroup of a wreath product from any linear
identity. This construction seems to be of interest on its own, since it
allows us to work with identities in a finite group instead of an (infinite)
free groupoid. Finally, we present the shortest known equational basis for
Boolean algebras, based on $\eqref{Eq:H}$.

\subsection{Related work}

The identity \eqref{Eq:H} was studied by Thedy \cite{Th1967} for rings. It
appears as identity (10) in \cite{JKVar}. Hossz\'u \cite{Ho} showed that a
quasigroup satisfying $(\ref{Eq:H})$ is an abelian group.

Kleinfeld \cite{Kl1978} investigated the left modular identity for rings.
Belousov \cite{Be} and Je\v{z}ek and Kepka \cite{JKVar, JKEq} worked with
linear identities in the variety of quasigroups. Equational theories of some
linear identities are studied in \cite{JKLin}. There is an extensive
bibliography of early papers on balanced identities (especially medial
groupoids) in the monograph \cite{JKMed}.

\section{Systematic Names for Linear Identities}\label{Sc:SN}

\noindent Although most proofs in this paper are easy to understand
intuitively, a systematic notation helps to write them down formally.

\subsection{Labelling bracketings}
Products of $n$ factors can be represented as labelled binary trees, or as
groupoid terms of length $n$. Unlabelled binary trees correspond to bracketings
of factors in a product. The \emph{length} of a bracketing is the number of
leaves in the corresponding tree.

Products of $n$ factors can be bracketed in $C_n$ ways, where $C_n$ is the
\emph{$n$th Catalan number} defined by the recurrence relation
\begin{equation}\label{Eq:CNR}
    C_1=1,\quad C_2=1,\quad
    C_n=C_1C_{n-1}+C_2C_{n-2}+\cdots +C_{n-2}C_2+C_{n-1}C_1,
\end{equation}
which is equivalent to the explicit formula
\begin{equation}\label{Eq:CNE}
    C_{n+1}=\frac{1}{n+1}\binom{2n}{n}.
\end{equation}
See \cite{vLW} for more on Catalan numbers and Table \ref{Tb:Count} for the
first few values $C_n$.

Given a bracketing $t$ of length $n$ it is therefore possible to assign a
unique name $b(t)$ to $t$ so that $0\le b(t)< C_n$. One way of doing this is
to (a) represent each bracketing by a sequence of symbols `(' (left
parenthesis), `)' (right parenthesis) and `$\circ$' (placeholder), (b)
introduce a total order on the three symbols, (c) extend this total order
lexicographically to a total order on all bracketings of given length.

In this paper, we will label bracketings as follows:

When $t$ is a bracketing of length $1$, let $b(t)=0$. When $t=t_\lambda
t_\rho$ is a bracketing of length $n>1$ that is a product of a bracketing
$t_\lambda$ of length $n-m$ and a bracketing $t_\rho$ of length $m$, let
\begin{equation}\label{Eq:Code}
    b(t) = \left(\sum_{i=1}^{m-1} C_iC_{n-i}\right) +
        b(t_\rho)C_{n-m} + b(t_\lambda).
\end{equation}
Thus, the function $b$ first counts all bracketings whose top two products
are of length $n-1$ and $1$, respectively, then moves on to all bracketings
whose top two products are of length $n-2$ and $2$, respectively, and so on.
To see that $b$ is a bijection, we prove that the bracketing $t$ can be
reconstructed from $b(t)$.

When $n=1$, $t$ is determined by $b(t)=0$. Assume that $b(s)$ determines $s$
uniquely for all bracketings $s$ of length less than $n$. Let $m$ be the
biggest integer such that $\sum_{i=1}^{m-1}C_{i}C_{n-i}\le b(t)$. Then $m$ is
the length of $t_\rho$, and $d = b(t_\rho)C_{n-m} + b(t_\lambda)$ is
therefore known. Since $b(t_\lambda)< C_{n-m}$ by the induction hypothesis,
we have $b(t_\rho) = \lfloor d/C_{n-m}\rfloor$,
$b(t_\lambda)=d-b(t_\rho)C_{n-m}$, thus reconstructing $t_\lambda$ and
$t_\rho$ from $b(t)$.

\begin{example} Here are the first $8$ of the $C_5=14$ bracketings of length
$5$:
\begin{align*}
    &(((\circ\circ)\circ)\circ)\circ = 0,\quad
    ((\circ(\circ\circ))\circ)\circ = 1,\quad
    ((\circ\circ)(\circ\circ))\circ = 2,\quad
    (\circ((\circ\circ)\circ))\circ = 3,\\
    &(\circ(\circ(\circ\circ)))\circ = 4,\quad
    ((\circ\circ)\circ)(\circ\circ) = 5,\quad
    (\circ(\circ\circ))(\circ\circ)=6,\quad
    (\circ\circ)((\circ\circ)\circ)=7.
\end{align*}
Note that the labelling does not agree with the lexicographic order.
\end{example}

\subsection{Naming linear groupoid identities} Let $u=v$ be a
linear groupoid identity of length $2n$. Let $b(u)$ be the label of the
bracketing of $u$, $b(v)$ the label of the bracketing of $v$, and $f\in S_n$
the permutation that must be applied to the variables of $u$ so that they
become ordered as in $v$. Since every variable occurs exactly once on both
sides, the permutation $f$ is uniquely determined. We can hence identify the
identity $u=v$ with the quadruple $(n$, $b(u)$, $b(v)$, $f)$, which we call
the \emph{name} of $u=v$. In order to save space, we write
$\idname{n}{b(u)}{b(v)}{f}$ instead of $(n$, $b(u)$, $b(v)$, $f)$, or even
$\sidname{b(u)}{b(v)}{f}$, when $n$ is clear from the context. The variety of
groupoids defined by a single linear identity $\sidname{i}{j}{f}$ will also
be denoted by $\sidname{i}{j}{f}$.

\begin{example}
The identity $((xy)u)v=(xu)(vy)$ has name $\idname{4}{0}{2}{(2,4,3)}$.
\end{example}

\begin{remark} The notation can be extended to arbitrary balanced groupoid
identities. However, the permutation $f$ is then not necessarily uniquely
determined. It can be assigned canonically by imposing a total order on
permutations.
\end{remark}

\subsection{Counting linear identities}

Upon interchanging the left hand side and the right hand side of a linear
identity $\idname{n}{i}{j}{f}$, we obtain the identity
$\idname{n}{j}{i}{f^{-1}}$. Naturally, we consider these two identities to be
the same.

We call a linear identity \emph{trivial} if it is of the form
$\sidname{i}{i}{()}$, where $()$ is the identity permutation.

For $n$, $m>0$, let $s_{n,m}$ denote the number of elements of order $m$ in
the symmetric group $S_n$.

\begin{lemma}\label{Lm:Count}
There are
\begin{equation}\label{Eq:Count}
    \frac{C_n}{2}\left( C_nn! + 1 + s_{n,2}\right)
\end{equation}
linear identities of length $2n$. There are
\begin{equation}\label{Eq:NTCount}
    \frac{C_n}{2}\left( C_nn! - 1 + s_{n,2}\right)
\end{equation}
nontrivial linear identities of length $2n$.
\end{lemma}
\begin{proof}
In order to construct a linear identity $\sidname{i}{j}{f}$ of length $2n$,
we can choose each of the bracketings $i$, $j$ in $C_n$ ways, and the
permutation $f$ in $n!$ ways. We do not distinguish between
$\sidname{i}{j}{f}$ and $\sidname{j}{i}{f^{-1}}$; hence the factor $1/2$.
However, before we divide by $1/2$, we must add all identities
$\sidname{i}{j}{f}$ for which $\sidname{i}{j}{f}$ and
$\sidname{j}{i}{f^{-1}}$ are the same. This happens if and only if $i=j$ and
$f^2=1$. Since there are $1+s_{n,2}$ permutations $f$ of $S_n$ with $f^2=1$,
we have proved \eqref{Eq:Count}. Equation \eqref{Eq:NTCount} follows from
\eqref{Eq:Count} upon subtracting the $C_n$ trivial identities
$\sidname{i}{i}{()}$.
\end{proof}

For the sake of completeness, we give the formula for $s_{n,2}$, which is
certainly well known.

\begin{lemma}
There are
\begin{displaymath}
    s_{n,2}=\sum_{1\le
    2m\le n} \binom{n}{2m}\cdot 1\cdot 3\cdot 5\cdots (2m-1)
\end{displaymath}
involutions in $S_n$.
\end{lemma}
\begin{proof}
All involutions of $S_n$ can be obtained as follows: Select an even number
$0<2m\le n$ of elements. Split the $2m$ elements into $m$ pairs, each
corresponding to some transposition $(a,b)$. An easy induction shows that the
number of ways in which $2m$ elements can be split into pairs (equivalently,
the number of $1$-factorizations of the complete graph on $2m$ vertices) is
$1\cdot 3\cdot 5\cdots (2m-1)$.
\end{proof}

\section{Identity-hedrons and implications among linear identities}

\subsection{Free groupoids}

\noindent The \emph{absolutely free groupoid} $A_n$ on generators $x_1$,
$\dots$, $x_n$ consists of all groupoid terms formed from $x_1$, $\dots$,
$x_n$, i.e., of all words $u=x_{i_1}x_{i_2}\cdots x_{i_m}$ bracketed in some
way, where $m\ge 0$, $i_j\in\{1,\dots,n\}$ for $1\le j\le m$. The product of
two terms $u$, $v\in A_n$ is the term $uv$.

Let $\varphi$ be a groupoid identity. Define a binary relation $\sim$ on
$A_n$ by $u\sim v$ if and only if $v$ is obtained from $u$ by a single
application of the identity $\varphi$. Let $\equiv$ be the reflexive and
transitive closure of $\sim$ on $A_n$. Then $\equiv$ is a congruence of
$A_n$, and $F=A_n/{\equiv}$ (also denoted by $A_n/\varphi$) is the \emph{free
groupoid} with $n$ generators satisfying $\varphi$. Elements of $F$
(equivalence classes of $A_n$) will be denoted by $[u]$, where $u\in A_n$.

Let $\varphi$, $\psi$ be two linear groupoid identities. We say that
$\varphi$ \emph{implies} $\psi$ if every groupoid satisfying $\varphi$ also
satisfies $\psi$.

\begin{theorem}\label{Th:Free}
Let $\varphi$, $\psi$ be two linear groupoid identities. Assume that $\psi$
is the identity $u=v$ and that it is of length $2n$. Then $\varphi$ implies
$\psi$ if and only if $[u]=[v]$ in $A_n/\varphi$.
\end{theorem}
\begin{proof} If $[u]=[v]$ in $A_n/\varphi$ then $\psi$ is obtained by a
repeated application of $\varphi$, and hence every groupoid satisfying
$\varphi$ also satisfies $\psi$.

Assume that $[u]\ne[v]$ in $A_n/\varphi$. Then $A_n/\varphi$ is a groupoid
satisfying $\varphi$ but not $\psi$, and hence $\varphi$ does not imply
$\psi$.
\end{proof}

\subsection{Identity-hedrons}

\noindent When viewed as a transformation of terms in an absolutely free
groupoid, the primary effect of a linear identity is to change the bracketing
of a given product and, at the same time, to permute the factors. This leads
us to the notion of an \emph{identity-hedron}, that we introduce by means of
an example.

Consider the linear identity \eqref{Eq:H}. Let $X$ be the set of all
bracketings of length $4$. As in Section \ref{Sc:SN}, we can identify $X$
with the set $\{0,\dots,4\}$, since $C_4=5$.

\begin{figure}[ht]\begin{center}\input{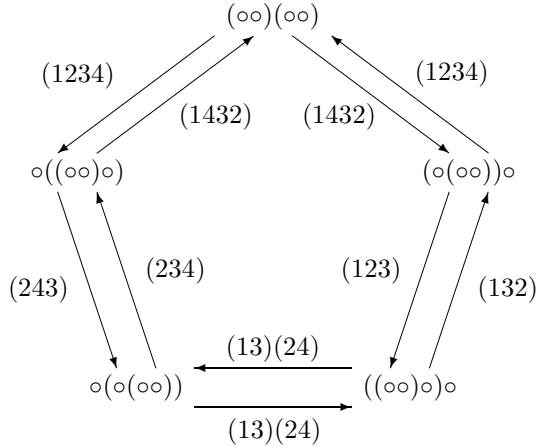}\end{center}
\caption{The identity $(xy)z=y(zx)$ and terms of length $4$}\label{Fg:New}
\end{figure}

Let $u=u_1u_2u_3u_4$ be a word bracketed in some way. Then the rule
\eqref{Eq:H} can be applied to it in several ways to yield another term. For
instance, when $u=((u_1u_2)u_3)u_4$, we can apply \eqref{Eq:H} in two ways to
obtain the terms $(u_2(u_3u_1))u_4$ and $u_3(u_4(u_1u_2))$, respectively.
Thus every application of \eqref{Eq:H} to a term $u$ is fully described by
the change in bracketing of $u$ and by the permutation of the letters $u_1$,
$\dots$, $u_4$. We can represent any such application of \eqref{Eq:H} by a
labelled arrow. Upon collecting all such arrows, we obtain an
\emph{identity-hedron}, as in Figure \ref{Fg:New}. (Here, our terminology is
analogous to \emph{associahedrons}. See \cite{Stasheff}.)

\setlength{\unitlength}{1.3mm}

Note that the information in an identity-hedron is redundant, since the two
arrows pointing in opposite directions are labelled by permutations that are
inverse to each other. We will exploit this redundancy later.

Naturally, the identity \eqref{Eq:H} can be applied more than once. This
corresponds to a journey through Figure $\ref{Fg:New}$ along a path of
arrows. Upon completing the journey, we are left with a permutation (obtained
by composing the permutations along the arrows), and hence with some linear
identity determined by the starting bracketing of the path, the terminating
bracketing of the path, and by the above permutation.

For instance, starting at bracketing $0=((\circ\circ)\circ)\circ$ and
travelling counterclockwise, we see that $((ab)c)d)=(b(ca))d = (db)(ca) =
a((db)c) = a(b(cd)) = ((cd)a)b$. We have returned to the same bracketing but
the order of the factors is different. We could have calculated the order of
the factors directly by rearranging $abcd$ according to the permutation
$(13)(24)(243)(1234)(1234)(132)=(13)(24)$. The linear identity corresponding
to this journey is thus $\idname{4}{0}{0}{(13)(24)}$.

It should now be clear how to construct an identity-hedron for any balanced
identity $\varphi$ and any length of terms $m$. We will denote the
corresponding identity-hedron by $H(\varphi, m)$.

We now make the anticipated connection between implications and journeys
through identity-hedrons.

\begin{theorem}\label{Th:Hedron}
Let $\varphi$, $\psi$ be linear identities.
Then $\varphi$ implies $\psi$ if and only if $\psi$ corresponds to a journey
through the identity-hedron $H(\varphi,m)$, where $2m$ is the length of
$\psi$.
\end{theorem}
\begin{proof}
Let $u$, $v$ be two terms in the absolutely free groupoid $A_m$ such that
$u=v$ is $\psi$. Then $[u]=[v]$ holds in $A_m/\varphi$ if and only if there
is a journey through $H(\varphi,m)$ that yields $\psi$. The rest follows from
Theorem \ref{Th:Free}.
\end{proof}

\section{Inclusions between Varieties of Groupoids Defined by a Linear
Identity}\label{Sc:Class}

\begin{theorem}\label{Th:Class}
 Let $n\ge 3$ be an integer and let $\sidname{i}{j}{f}$,
$\sidname{r}{s}{g}$ be two distinct, nontrivial linear identities of length
$2n$. Then $\sidname{i}{j}{f}$ implies $\sidname{r}{s}{g}$ if and only if
$i=j=r=s$ and $g=f^k$ for some $k\ne 0$.
\end{theorem}
\begin{proof}
By Theorem \ref{Th:Hedron}, $\sidname{i}{j}{f}$ implies $\sidname{r}{s}{g}$
if and only if $\sidname{r}{s}{g}$ is the result of a journey in
$H=H(\sidname{i}{j}{f},n)$. Without loss of generality, let $i\le j$, $r\le
s$.

Since the two identities in question are of the same length, the
identity-hedron $H$ is easy to describe: it consists of several isolated
bracketings and one pair of mutually inverse arrows connecting bracketings
$i$ and $j$ (when $i=j$, the arrows are loops).

When $i\ne j$, any journey through $H$ is of the form: (a)
$\sidname{i}{j}{f}$, (b) $\sidname{j}{i}{f^{-1}}$, (c) $\sidname{i}{i}{()}$,
or $\sidname{j}{j}{()}$. None of these identities is $\sidname{r}{s}{g}$
since: (a) $\sidname{i}{j}{f}\ne\sidname{r}{s}{g}$, (b) $i<j$, $r\le s$, (c)
$g$ is nontrivial.

When $i=j$ then any journey through $H$ is of the form $\sidname{i}{i}{f^k}$,
and we are done.
\end{proof}

\begin{example}
By Theorem \ref{Th:Class}, the identity
\begin{displaymath}
    \idname{4}{3}{3}{(1,2,3,4)}\ =\ ``x((yz)u)=u((xy)z)"
\end{displaymath}
implies the identity
\begin{displaymath}
    \idname{4}{3}{3}{(1,3)(2,4)}\ =\ ``x((yz)u)=z((ux)y)",
\end{displaymath}
but the two identities are not equivalent. This is also witnessed by the
groupoid with the following multiplication table:
\begin{displaymath}
\begin{array}{c|ccc}
    &0&1&2\\
    \hline
    0&0&1&2\\
    1&2&0&1\\
    2&1&2&0
\end{array}
\end{displaymath}
\end{example}

We can now easily count the varieties of groupoids defined by a single linear
identity of length $2n$. Recall that $s_{n,m}$ denotes the number of elements
of order $m$ in the group $S_n$. Let $\varphi$ be the \emph{Euler function},
i.e., $\varphi(m)$ is the number of positive integers less than $m$ that are
relatively prime to $m$.

\begin{proposition}\label{Pr:VarCount}
There are
\begin{equation}\label{Eq:VarCount}
    \mathcal L(n) = 1 + \binom{C_n}{2}n!
    + C_n\sum_{m\ge2} \frac{s_{n,m}}{\varphi(m)}
\end{equation}
varieties of groupoids defined by a single linear identity of length $2n$.
\end{proposition}
\begin{proof}
There is $1$ trivial variety (all groupoids). The second summand of
\eqref{Eq:VarCount} accounts for all linear identities $\sidname{i}{j}{f}$
with $i<j$. It remains to count the varieties defined by some nontrivial
$\sidname{i}{i}{f}$. We know from Theorem \ref{Th:Class} that
$\sidname{i}{i}{f}=\sidname{j}{j}{g}$ (as varieties) if and only if $i=j$,
$f=g^k$ and $g=f^l$ for some $k\ne 0\ne l$. There are $C_n$ bracketings $i$
of length $n$. If $f\in S_n$ is a permutation of order $m$, it gives rise to
a cyclic subgroup $G\le S_n$ of order $m$. The permutations $g$ satisfying
$g=f^k$, $f=g^l$ for some $k\ne 0\ne l$ are then precisely the generators of
$G$. It is well known that a cyclic group of order $m$ has $\varphi(m)$
generators.
\end{proof}

\begin{remark} The summand $s_{n,m}/\varphi(m)$ counts the number of cyclic
subgroups of order $m$ in $S_n$, and, therefore, the sum $\sum_{m\ge 2}
\frac{s_{n,m}}{\varphi(m)}$ in $(\ref{Eq:VarCount})$ is the number of
nontrivial cyclic subgroups of $S_n$.
\end{remark}

\begin{table}
\caption{The number $\mathcal L(n)$ of varieties of groupoids defined by a
single linear identity of length $2\le n\le 6$, together with all constants
needed to evaluate $\mathcal L(n)$, based on Proposition
\ref{Pr:VarCount}.}\label{Tb:Count}
\begin{displaymath}
    \begin{array}{r|rrrrr}
    s_{n,m}& 2&3&4&5&6\\
    \hline
    2&1& & & &\\
    3&3&2 & & & \\
    4&9&8&6& & \\
    5&25&20&30&24&20\\
    6&75&80&180&144&240
    \end{array}
    \quad\quad
    \begin{array}{r|rrrrr}
    n&2&3&4&5&6\\
    \hline
    \varphi(n)&1&2&2&4&2\\
    C_n&1&2&5&14&42\\
    \mathcal L(n)&2&15&321&11,845&635,083
    \end{array}
\end{displaymath}
\end{table}

Table \ref{Tb:Count} gives the values of $\mathcal L(n)$, for $2\le n\le 6$.

\begin{theorem}\label{Th:nm} Assume that $f\in S_n$, $g\in S_m$ are
nonidentity permutations, $n\ne m$. Then the varieties $\idname{n}{i}{j}{f}$,
$\idname{m}{r}{s}{g}$ are not the same.
\end{theorem}
\begin{proof} Without loss of generality, let $n<m$. Then the
identity-hedron $H(\idname{m}{r}{s}{g},n)$ contains no arrows, and hence
$\idname{m}{r}{s}{g}$ does not imply $\idname{n}{i}{j}{f}$, by Theorem
\ref{Th:Hedron}.
\end{proof}

\subsection{The fourteen varieties of length $6$}\label{Sc:14}

\noindent Note that $14$ is both the number of nontrivial varieties defined by
a single linear identity of length $6$ (Table \ref{Tb:Count}), and the number
of nontrivial linear identities of length $6$ (Lemma \ref{Lm:Count}). For the
convenience of the reader, these $14$ identities can be found in Table
\ref{Tb:Varieties}.

\begin{table}
\caption{The $14$ nontrivial varieties of groupoids defined by a single linear
identity of length $6$. We omit ``3'' from their systematic
names.}\label{Tb:Varieties}
\begin{displaymath}
\begin{array}{cclc}
    \text{identity}&\text{systematic name}&\text{is equivalent to}
    &\text{remark}\\
    \hline
    (xy)z=(yx)z&\sidname{0}{0}{(1,2)}&&\\
    (xy)z=(zy)x&\sidname{0}{0}{(1,3)}&&\text{right modular groupoids}\\
    (xy)z=(xz)y&\sidname{0}{0}{(2,3)}&&\text{right permutable groupoids}\\
    (xy)z=(zx)y&\sidname{0}{0}{(1,2,3)}&\sidname{0}{0}{(1,3,2)}&\\
    (xy)z=x(yz)&\sidname{0}{1}{()}&\sidname{1}{0}{()}&\text{semigroups}\\
    (xy)z=y(xz)&\sidname{0}{1}{(1,2)}&\sidname{1}{0}{(1,2)}&\\
    (xy)z=z(yx)&\sidname{0}{1}{(1,3)}&\sidname{1}{0}{(1,3)}&\\
    (xy)z=x(zy)&\sidname{0}{1}{(2,3)}&\sidname{1}{0}{(2,3)}&\\
    (xy)z=z(xy)&\sidname{0}{1}{(1,2,3)}&\sidname{1}{0}{(1,3,2)}&\\
    (xy)z=y(zx)&\sidname{0}{1}{(1,3,2)}&\sidname{1}{0}{(1,2,3)}&
        \text{Eq. \eqref{Eq:H}}\\
    x(yz)=y(xz)&\sidname{1}{1}{(1,2)}&&\\
    x(yz)=z(yx)&\sidname{1}{1}{(1,3)}&&\text{left modular groupoids}\\
    x(yz)=x(zy)&\sidname{1}{1}{(2,3)}&&\text{left permutable groupoids}\\
    x(yz)=z(xy)&\sidname{1}{1}{(1,2,3)}&\sidname{1}{1}{(1,3,2)}&
\end{array}
\end{displaymath}
\end{table}

For each identity $\sidname{i}{j}{f}$ of Table \ref{Tb:Varieties} we now
construct a finite groupoid satisfying $\sidname{i}{j}{f}$ but not any other of
the remaining $13$ identities. The multiplication tables of these groupoids are
gathered in Figure \ref{Fg:14}. The multiplication table of a groupoid
satisfying $\sidname{i}{j}{f}$ is labelled by $\sidname{i}{j}{f}$. All $m\times
m$ multiplication tables of Figure \ref{Fg:14} have rows and columns labelled
by $0$, $\dots$, $m-1$, in this order.

\begin{figure}
\setlength\arraycolsep{3pt}
\begin{footnotesize}
\begin{gather*}
    \begin{array}{c}
        \sidname{1}{1}{(1,2)}\\
        \begin{array}{|ccc}
            \hline
            0&1&0\\
            0&1&2\\
            1&1&1
        \end{array}
    \end{array}
\\
\begin{array}{cccccc}
    \begin{array}{c}
        \sidname{0}{0}{(1,2)}\\
        \begin{array}{|cccc}
            \hline
            0&0&3&2\\
            0&0&3&0\\
            0&0&3&2\\
            0&0&3&2
        \end{array}
    \end{array}
    &
    \begin{array}{c}
        \sidname{0}{0}{(2,3)}\\
        \begin{array}{|cccc}
            \hline
            0&3&3&3\\
            1&1&1&1\\
            2&0&0&0\\
            3&2&2&2
        \end{array}
    \end{array}
    &
    \begin{array}{c}
        \sidname{0}{1}{()}\\
        \begin{array}{|cccc}
            \hline
            0&0&2&2\\
            1&1&3&3\\
            0&0&2&2\\
            1&1&3&3
        \end{array}
    \end{array}
    &
    \begin{array}{c}
        \sidname{0}{1}{(1,2)}\\
        \begin{array}{|cccc}
            \hline
            0&1&2&2\\
            0&1&3&3\\
            0&1&2&2\\
            0&1&2&2
        \end{array}
    \end{array}
    &
    \begin{array}{c}
        \sidname{0}{1}{(2,3)}\\
        \begin{array}{|cccc}
            \hline
            0&0&3&0\\
            1&1&1&1\\
            2&2&2&2\\
            0&0&3&0
        \end{array}
    \end{array}
    &
    \begin{array}{c}
        \sidname{1}{1}{(2,3)}\\
        \begin{array}{|cccc}
            \hline
            0&0&3&0\\
            1&1&0&1\\
            1&1&0&1\\
            0&0&0&0
        \end{array}
    \end{array}
\end{array}
\\
\begin{array}{cc}
    \begin{array}{c}
        \sidname{0}{0}{(1,3)}\\
        \begin{array}{|ccccc}
            \hline
            0&3&4&1&2\\
            2&1&0&4&3\\
            3&4&2&0&1\\
            4&2&1&3&0\\
            1&0&3&2&4
        \end{array}
    \end{array}
    &
    \begin{array}{c}
        \sidname{1}{1}{(1,3)}\\
        \begin{array}{|ccccc}
            \hline
            0&3&0&1&1\\
            1&0&1&3&3\\
            0&3&0&1&1\\
            3&1&4&0&0\\
            3&1&3&0&0
        \end{array}
    \end{array}
\end{array}
\\
\begin{array}{cc}
    \begin{array}{c}
        \sidname{0}{1}{(1,3)}\\
        \begin{array}{|cccccc}
            \hline
            3&0&1&3&3&4\\
            4&3&1&3&0&3\\
            5&5&3&3&1&1\\
            3&3&3&3&3&3\\
            3&4&5&3&3&0\\
            0&3&5&3&4&3
        \end{array}
    \end{array}
    &
    \begin{array}{c}
        \sidname{0}{1}{(1,2,3)}\\
        \begin{array}{|cccccc}
            \hline
            3&3&3&3&4&3\\
            4&3&3&4&3&3\\
            3&4&3&5&4&3\\
            3&4&5&3&3&3\\
            4&3&4&3&3&3\\
            3&3&3&3&3&3
        \end{array}
    \end{array}
\end{array}
\\
\begin{array}{ccc}
    \begin{array}{c}
        \sidname{0}{0}{(1,2,3)}\\
        \begin{array}{|ccccccccc}
            \hline
            2&3&4&4&4&4&8&4&4\\
            2&3&6&6&4&4&4&4&4\\
            5&5&4&4&4&4&4&4&4\\
            5&5&7&7&4&4&4&4&4\\
            4&4&4&4&4&4&4&4&4\\
            4&7&4&4&4&4&4&4&4\\
            7&7&4&4&4&4&4&4&4\\
            4&4&4&4&4&4&4&4&4\\
            4&4&4&4&4&4&4&4&4
        \end{array}
    \end{array}
    &
    \begin{array}{c}
        \sidname{0}{1}{(1,3,2)}\\
        \begin{array}{|ccccccccc}
            \hline
            3&3&3&4&4&7&8&4&4\\
            4&4&6&4&4&7&4&4&4\\
            5&3&5&4&4&7&8&4&4\\
            4&4&7&4&4&4&4&4&4\\
            4&4&4&4&4&4&4&4&4\\
            4&8&7&4&4&4&4&4&4\\
            4&4&4&4&4&4&4&4&4\\
            4&4&4&4&4&4&4&4&4\\
            4&4&4&4&4&4&4&4&4
        \end{array}
    \end{array}
    &
    \begin{array}{c}
        \sidname{1}{1}{(1,2,3)}\\
        \begin{array}{|ccccccccc}
            \hline
            1&3&1&4&4&3&4&4&7\\
            4&4&8&4&4&7&4&4&4\\
            5&3&5&7&4&3&4&4&7\\
            4&4&6&4&4&4&4&4&4\\
            4&4&4&4&4&4&4&4&4\\
            4&4&8&4&4&7&4&4&4\\
            4&4&4&4&4&4&4&4&4\\
            4&4&4&4&4&4&4&4&4\\
            4&4&4&4&4&4&4&4&4
        \end{array}
    \end{array}
\end{array}
\end{gather*}
\end{footnotesize}
\caption{The groupoid man. Smallest distinguishing examples for the $14$
nontrivial varieties of groupoids defined by a single linear identity of length
$6$.}\label{Fg:14}
\end{figure}

All examples in Figure \ref{Fg:14} are as small as possible. They were found by
Mace 4 \cite{Mc}. The groupoid $\sidname{0}{0}{(1,2,3)}$ was hardest to find;
it took Mace 4 about $5$ hours on a Pentium 3 machine with 765 megabytes of
RAM.

\begin{corollary} The $14$ nontrivial varieties defined by a single linear
identity of length $6$ are in a general position with respect to inclusion,
i.e., none of these varieties is contained in the union of the remaining $13$
varieties.
\end{corollary}

\begin{remark} Je\v{z}ek and Kepka determined that there are $11$ varieties of
quasigroups defined by a single linear identity of length $\le 6$, and found
all inclusions among them \cite[Theorem 1.8]{JKVar}. Kirnasovsky \cite{Ki}
studied the same problem for length $\le 8$.
\end{remark}

\subsection{Linear identities of length $8$ implied by $\eqref{Eq:H}$}\label{Sc:4}

\noindent Let us answer a question posed in \cite{HJM1993}.

\begin{proposition} Exactly $45$ out of the $320$ nontrivial
linear identities of length $8$ are implied by $\eqref{Eq:H}$. These identities
can be found with the aid of Figure $\ref{Fg:New}$.
\end{proposition}
\begin{proof}
Recall that upon completing one counterclockwise round in Figure \ref{Fg:New}
starting at bracketing $0$, the $4$ symbols are permuted according to
$(13)(24)$. We claim that all permutations $g$ corresponding to one
counterclockwise round in Figure \ref{Fg:New} are involutions. This is easy
to see, since the permutation obtained by starting at bracketing $i$ is a
conjugate of the permutation obtained by starting at bracketing $i-1$ (cf.,
for $i=1$, we get $(132)(13)(24)(132)^{-1}$).

We can now describe all journeys through Figure \ref{Fg:New}: Select two
bracketings $i$, $j$ out of the $5$ possible bracketings, allowing $i=j$. The
shortest nonempty counterclockwise path from $i$ to $j$ yields some identity
$\sidname{i}{j}{f}$. Upon extending this path by another complete
counterclockwise round, we obtain identity $\sidname{i}{j}{gf}$ that is
different from $\sidname{i}{j}{f}$ (since $g$ is an involution). The two
bracketings $i$, $j$ can be chosen in $5\cdot 5 = 25$ ways. Hence we obtain
$50$ identities following the above procedure. Five of these identities are
trivial (those corresponding to two full rounds).

We are done by Theorem \ref{Th:Hedron}.
\end{proof}

We verified by the finite model builder Mace 4 \cite{Mc} that the groupoid
$A\times B$ defined by
\begin{displaymath}
\begin{array}{c|ccccccccc}
  A&0&1&2&3&4&5&6&7&8\\
  \hline
  0&2&3&5&7&5&6&8&8&8\\
  1&4&3&6&6&6&6&8&8&8\\
  2&5&5&6&8&8&8&8&8&8\\
  3&6&6&6&8&8&8&8&8&8\\
  4&7&7&6&8&8&8&8&8&8\\
  5&6&8&8&8&8&8&8&8&8\\
  6&8&8&8&8&8&8&8&8&8\\
  7&6&8&8&8&8&8&8&8&8\\
  8&8&8&8&8&8&8&8&8&8
\end{array}\quad\quad
\begin{array}{c|cccccccc}
    B&0&1&2&3&4&5&6&7\\
    \hline
    0&2&3&2&7&2&6&2&2\\
    1&4&2&2&2&5&2&2&6\\
    2&2&2&2&2&2&2&2&2\\
    3&2&5&2&2&6&2&2&2\\
    4&7&2&2&6&2&2&2&2\\
    5&6&2&2&2&2&2&2&2\\
    6&2&2&2&2&2&2&2&2\\
    7&2&6&2&2&2&2&2&2
\end{array}
\end{displaymath}
satisfies \eqref{Eq:H} but none of the remaining $320-45$ linear identities
of length $8$ not implied by \eqref{Eq:H}.

\section{Ultimately AC-nice groupoids and labelled binary trees}\label{Sc:UNG}

\noindent As in \cite{St} and \cite{Mu}, two groupoid terms are said to be
\emph{AC-identical} if one is obtained from the other by a repeated
application of associativity and commutativity.

For an integer $m>1$, we then say that a groupoid $G$ is \emph{$m$AC-nice} if
any two products of the same $m$ elements of $G$ yield the same element of
$G$. $m$AC-nice groupoids are called \emph{$m$-nice} in \cite{HJM1993}.

By \cite[Lemma 2.2]{HJM1993}, every groupoid satisfying \eqref{Eq:H} is
$5$AC-nice. By \cite[Lemma 2.3]{HJM1993}, an $m$AC-nice groupoid is
$(m+1)$AC-nice, provided $m\ge 3$. Note that $2$AC-nice groupoids are
precisely commutative groupoids, and $3$AC-nice groupoids are groupoids that
are commutative and associative. $3$AC-niceness therefore does not follow
from $2$AC-niceness.

It thus makes sense to say:

\begin{definition} A groupoid
$G$ is \emph{ultimately AC-nice} if it is $m$AC-nice for some $m\ge 3$. A
linear identity $\sidname{i}{j}{f}$ is \emph{ultimately AC-nice} if every
groupoid satisfying $\sidname{i}{j}{f}$ is ultimately AC-nice.
\end{definition}

In the last paragraph of \cite{HJM1993}, the authors of \cite{HJM1993} claim,
without proof, that there is a groupoid that satisfies $x(yz)=z(yx)$ but that
is not $5$AC-nice. We prove a general result (Theorem \ref{Th:Nice}) along
similar lines: \emph{the only ultimately AC-nice linear identity is
$(\ref{Eq:H})$.}

\setlength{\unitlength}{1.0mm}
\begin{figure}[ht]\begin{center}\input{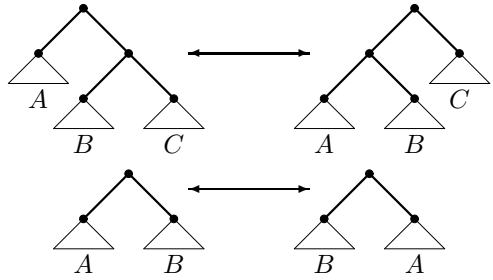}\end{center}
\caption{Associativity and commutativity as transformations of labelled binary
trees}\label{Fg:AC}
\end{figure}

This result can be visualized in terms of transformations of labelled binary
trees as follows:

First notice that a linear identity $\sidname{i}{j}{f}$ is ultimately AC-nice
if and only if all free groupoids satisfying $\sidname{i}{j}{f}$ are
ultimately AC-nice. One application of a linear identity to a word in the
absolutely free groupoid can be depicted by two labelled binary trees. Figure
\ref{Fg:AC} shows this for the associative law and for the commutative law.
Since two groupoid products with the same factors coincide in the presence of
associativity and commutativity, we see that given two labelled binary trees
$T_1$, $T_2$ with the same $n$ leaves, it is possible to obtain $T_2$ from
$T_1$ by finitely many applications of the two laws. Can the same feat be
achieved by a single linear identity $\varphi$, at least for sufficiently
large trees? This is precisely the question whether $\varphi$ is ultimately
AC-nice, and we answer it in Theorem \ref{Th:Nice}. The proof of Theorem
\ref{Th:Nice} is split into several steps:

\begin{lemma}\label{Lm:FixEnd}
Let $\idname{n}{i}{j}{f}$ be a linear identity such that $f\in S_n$ satisfies
$f(1)=1$ or $f(n)=n$. Then the free groupoid on two generators satisfying
$\idname{n}{i}{j}{f}$ is not ultimately AC-nice.
\end{lemma}
\begin{proof}
Let $A$ be the absolutely free groupoid on generators $x$, $y$. For $m\ge 3$,
consider the words $u=u_1\dots u_m$, $v=v_1\dots v_m\in A$ such that
$u_1=v_m=x$, $v_1=u_m=y$, $u_k=v_k=x$ for $1<k<m$. Assume that $f(1)=1$. No
matter what the bracketing of $u$ is, we see that no application of the
identity $\sidname{i}{j}{f}$ can move $u_1$ from the left-most position.
Since $u_1\ne v_1$, the products $u$, $v$ do not coincide in
$A/\sidname{i}{j}{f}$. Similarly when $f(n)=n$.
\end{proof}

\begin{lemma}\label{Lm:Comm}
The free commutative groupoid on one generator is not ultimately AC-nice.
\end{lemma}
\begin{proof}
Let $F$ be the free commutative groupoid with generator $x$. Define powers
$x^n$ by $x^1=x$, $x^n=xx^{n-1}$. Then for any even $m\ge 4$ we have $x^m\ne
(xx)^{m/2}$, since commutativity is not strong enough to split any of the
factors $xx$.
\end{proof}

\setlength{\unitlength}{1.0mm}
\begin{figure}[ht]\begin{center}\input{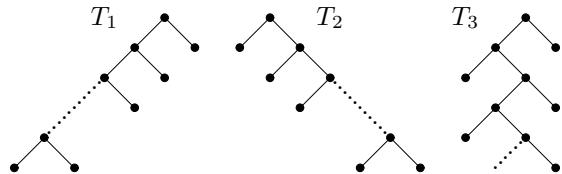}\end{center}
\caption{The proof of Theorem \ref{Th:Nice}}\label{Fg:3Trees}
\end{figure}

\begin{proposition}\label{Pr:LR}
Let $\idname{n}{i}{j}{f}$ be an ultimately AC-nice linear identity with $i\le
j$. Then $i=0$, $j=C_n-1$.
\end{proposition}
\begin{proof}
Let $F$ be the free groupoid on $1$ generator satisfying
$\idname{n}{i}{j}{f}$. Since $\idname{n}{i}{j}{f}$ is ultimately AC-nice, it
must be possible to transform the tree $T_1$ of Figure \ref{Fg:3Trees} into
the tree $T_2$ of the same Figure by a repeated application of
$\idname{n}{i}{j}{f}$, provided the two trees have the same number of leaves
and are sufficiently large. Note that $\sidname{i}{j}{f}$ is applicable to
$T_1$ if and only if the bracketing $i$ is of the form
$(\cdots((\circ\circ)\circ)\cdots)\circ$, i.e., if and only if $i=0$.
Similarly, $\sidname{i}{j}{f}$ is applicable to $T_2$ if and only if
$j=C_n-1$.
\end{proof}

\begin{proposition}\label{Pr:6}
The only ultimately AC-nice linear identity of length $6$ is $(\ref{Eq:H})$.
\end{proposition}
\begin{proof}
Proposition \ref{Pr:LR} eliminates $8$ identities of the form
$\sidname{i}{i}{f}$ from Table \ref{Tb:Varieties}. The identities
$\sidname{0}{1}{()}$, $\sidname{0}{1}{(1,2)}$, $\sidname{0}{1}{(2,3)}$ fix
either $1$ or $3$, and are therefore eliminated by Lemma \ref{Lm:FixEnd}.
Finally, the two identities $\sidname{0}{1}{(1,3)}$,
$\sidname{0}{1}{(1,2,3)}$ are consequences of the commutative law, and hence
are eliminated by Lemma \ref{Lm:Comm}.
\end{proof}

\begin{theorem}\label{Th:Nice} The only ultimately AC-nice linear identity is
$(\ref{Eq:H})$.
\end{theorem}
\begin{proof}
The only nontrivial linear identity of length $\le 4$ is the commutative law
$xy=yx$, which is not ultimately AC-nice by Lemma \ref{Lm:Comm}. Thanks to
Proposition \ref{Pr:6}, it suffices to consider ultimately AC-nice linear
identities of length $\ge 8$. Let $\idname{n}{i}{j}{f}$ be such an identity,
$i\le j$, $n\ge 4$. By Proposition \ref{Pr:LR}, we have $i=0$, $j=C_n-1$.
Consider the tree $T_3$ of Figure \ref{Fg:3Trees}. We claim that
$\idname{n}{i}{j}{f}$ is not applicable to $T_3$, no matter how large $T_3$
is. This is because it is impossible to make at least $n-1$ consecutive moves
to the left (or to the right) along the branches of $T_3$.
\end{proof}

\setlength{\unitlength}{1.0mm}
\begin{figure}[ht]\begin{center}\input{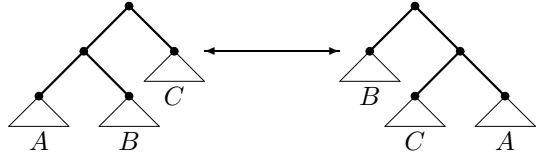}\end{center}
\caption{The identity $(\ref{Eq:H})$ as a transformation of labelled binary
trees}\label{Fg:HRule}
\end{figure}

The only ultimately AC-nice linear identity $(\ref{Eq:H})$ is visualized in
Figure \ref{Fg:HRule}.

Note that our proofs depend essentially on infinite (free) groupoids. Is this
dependence necessary?

\begin{conjecture} Let $\idname{n}{i}{j}{f}$ be a linear identity such that
every \emph{finite} groupoid satisfying $\idname{n}{i}{j}{f}$ is ultimately
AC-nice. Then $\idname{n}{i}{j}{f}$ is the identity $(\ref{Eq:H})$.
\end{conjecture}

\section{Wreath Products Associated with Linear Identities}\label{Sc:WP}

\noindent Let $\varphi$ be a linear identity, and $m>0$ an integer. By
composing arrows in the identity-hedron $H(\varphi,m)$, we can determine all
linear identities of length $2m$ implied by $\varphi$. Although it may seem
that it only makes sense to compose consecutive arrows of an identity-hedron,
we show below that it is possible to compose arbitrary arrows.

In this section, maps are applied to the right of their arguments, and
therefore composed from left to right.

\subsection{The associated wreath products}

Let us first recall some group-theoretical definitions:

Let $B$ be a group acting on another group $A$ via $a\mapsto a^b$, where
$a\in A$, $b\in B$. Then the \emph{semidirect product} $A\ltimes B$ is the
group defined on $A\times B$ by
$(a_1,b_1)(a_2,b_2)=(a_1a_2^{(b_1^{-1})},b_1b_2)$. We use $a_2^{(b_1^{-1})}$
rather than $a_2^{b_1}$ in the definition of a semidirect product because we
compose maps from left to right.

Let $B$ be a group acting on a set $X$, and let $A$ be another group. Then
$B$ also acts on the set $A^X$ of all maps from $X$ to $A$ via
$xf^b=x^{b^{-1}}f$, $f\in A^X$, $b\in B$, $x\in X$. The \emph{wreath product}
$\wreath{A}{B}$ of $A$ and $B$ is the semidirect product $A^X\ltimes B$ under
this action.

When $X$ is a finite set $\{1,\cdots,n\}$, the maps $A^X$ can be identified
with the direct product $A^n$, and the elements of $\wreath{A}{B}$ can be
represented as $((p_1,\dots,p_n),p)=((p_i),p)$, where $p_i\in A$, $p\in B$.
When $B$ is a subgroup of $S_n$ acting naturally on $X$, the multiplication in
$\wreath{A}{B}$ is described by the explicit formula
\begin{equation}\label{Eq:WP}
    ((a_i),a)\cdot((b_i),b) =((a_i)(b_i)^{(a^{-1})},ab)
     =((a_i b_{ia}),ab),
\end{equation}
where, in accordance with our conventions, $ia$ is the image of $i$ under $a$.

Let us return to linear identities.

Fix a linear identity $\varphi$. Let $m$ be a positive integer and let $X$ be
the set of all bracketings of length $m$, $X=\{0,\dots,C_m-1\}$. Let $B=S_X$,
$A=S_m$, and $W=\wreath{A}{B}$.

Consider the arrow leading from bracketing $i\in X$ to bracketing $j\in X$
labelled by $\pi\in A$ in the identity-hedron $H(\varphi,m)$. We will
represent this arrow and its inverse (exploiting the redundancy) by a single
element $((a_0,\dots,a_{C_m-1}),a)=((a_i),a)$ of $W$ by letting $a$ be the
transposition $(i,j)$, and by setting $a_i=\pi$, $a_j=\pi^{-1}$,
$a_r=id_{\{1,\dots,m\}}$ for $r\not\in\{i,j\}$.

We claim that the multiplication formula $\eqref{Eq:WP}$ then generalizes
composition of consecutive arrows (transformations). To see this, consider
the word $u$ bracketed according to $i$. Let $v$ be the word obtained from
$u$ when $((a_j),a)$ is applied to $u$. Since $u$ is bracketed according to
$i$, the permutations $a_j$, $j\ne i$ are irrelevant. Hence $v$ will be
bracketed according to $ia$ and the letters of $u$ will be reordered in $v$
according to $a_i$. Let $w$ be the word obtained from $v$ after $((b_j),b)$
is applied to $v$. Then $w$ is bracketed according to $iab$ and the letters
of $u$ will be reordered in $w$ according to $a_ib_{ia}$. This agrees with
\eqref{Eq:WP}.

\begin{definition} Given a linear identity $\varphi$ and a positive integer
$m$, let $W(\varphi,m)$ be the subgroup of $W$ generated by the elements
$((a_i),a)\in W$ corresponding to all arrows (and their inverses) in the
identity-hedron $H(\varphi,m)$, as described above.
\end{definition}

\subsection{Wreath products and AC-niceness}\label{Ss:WN}

\noindent We have managed to associate a certain subgroup $W(\varphi,m)$ of a
wreath product with a linear identity $\varphi$ and a positive integer $m$.
We now show how these wreath products can be used to express $m$AC-niceness
for $\eqref{Eq:H}$. Conceivably, $W(\varphi,m)$ will be useful in other
settings, too.

Let $G$ be a subgroup of $\wreath{A}{B}$. Then $G$ acts on $X$ via the
original action of $B$, i.e., $x^{((p_i),p)}=xp$, where $x\in X$, $p_i\in A$,
$p\in B$. For $x\in X=\{0$, $\dots$, $C_m-1\}$, let $G_x\le G$ be the
stabilizer of $x$ and $O_x\subseteq X$ the orbit of $x$ under this action of
$G$. Denote by $P_x$ the projection of $G_x$ onto the $x$th component of
$\wreath{A}{B}$, i.e., $P_x=\{p_x;\;((p_i),p)\in G_x\}\le A$.

\begin{proposition}\label{Pr:kN}
Let $\varphi$ be a linear identity and $m$ a positive integer. Then all
groupoids satisfying $\varphi$ are $m$AC-nice if and only if there is a
bracketing $x\in X$ such that $G=W(\varphi,m)$ satisfies $P_x=S_m$, $O_x=X$.
\end{proposition}
\begin{proof}
Let $x\in X$ be such that $P_x=A=S_m$ and $O_x=X$. Let $H$ be a groupoid
satisfying $\varphi$, and let $u$, $v$ be two products with the same $m$
factors. Since $O_x=X$, the bracketings of $u$, $v$ can be changed to $x$.
Let $u'$, $v'$ be the corresponding products bracketed according to $x$.
Since $P_x=A$, the factors of $u'$, $v'$ can be reordered freely without
changing the value of $u'$, $v'$. Hence $u=u'=v'=v$ and $H$ is $m$AC-nice.

For the converse, assume that every groupoid satisfying $\varphi$ is
$m$AC-nice. Let $F=A_m/\varphi$ be the free groupoid on $m$ generators
satisfying $\varphi$. By our assumption, $F$ is $m$AC-nice. Fix a bracketing
$x$. Consider any two words $u$, $v\in A_m$ bracketed according to $x$. Since
$F$ is $m$AC-nice, $[u]=[v]$ in $F$. In other words, $v$ can be obtained from
$u$ by a repeated application of $\varphi$. This shows $P_x=S_m$. We can show
similarly that $O_x=X$.
\end{proof}

\begin{example} Consider again Figure \ref{Fg:New}. Let $x$ be any of the $5$
bracketings. Clearly, $O_x=X$. Thanks to the shape of the transformation
diagram (a cycle), it is also easy to see that the only way to return to $x$
is to complete several clockwise or counterclockwise cycles around the
diagram. Let $\pi$ be the permutation of the four symbols obtained after one
counterclockwise round starting at $x$. Then $\pi^{-1}$ corresponds to one
clockwise round. Hence $P_x=\langle \pi\rangle$. Since $S_4$ is not cyclic,
we have $P_x\ne S_4$. Since $x$ was arbitrary, we have proved that
\eqref{Eq:H} is not $4$AC-nice, by Proposition \ref{Pr:kN}.
\end{example}

\subsection{Computing the associated wreath products in GAP}\label{Ss:GAP}

\noindent Are calculations in $W(\varphi,m)$ more convenient than those in
the free groupoid on $m$ generators satisfying $\varphi$? It depends.

The advantage of $W(\varphi,m)$ is that it is a finite group, and hence all
tools of computational group theory apply to it. Importantly, up to $C_m$
applications of $\varphi$ are encoded in a single element of $W(\varphi,m)$.
Also note that the elements of $W(\varphi,m)$ capture the essence of
$\varphi$; namely all possible applications of $\varphi$ to words of length
$m$, not the words themselves.

On the other hand, $W(\varphi,m)$ is huge. There are $m^m\cdot C_m$ terms of
length $m$ in the free groupoid on $m$ generators. In comparison, the size of
$W$ (of which $W(\varphi,m)$ is a subgroup) is $(m!)^{C_m}\cdot C_m!$,
eventually a much bigger number.

By Proposition \ref{Pr:kN}, $m$AC-niceness of a linear identity $\varphi$ can
be determined by the study of the (projections of) stabilizers and the orbits
of the action of $W(\varphi,m)$ on all bracketings of length $m$.

Since stabilizers and orbits of permutation groups are implemented
efficiently in GAP \cite{GAP}, we wrote a short library of functions that
verifies $m$AC-niceness for a given linear identity $\varphi$. The library is
available electronically \cite{OurLibrary}. We describe the main functions
here.

Given a positive integer $m$ and a linear identity $\varphi$, the function
\begin{displaymath}
    \text{\texttt{GeneratorsByIdentity(}}m,\ \varphi\text{\texttt{)}}
\end{displaymath}
returns the generators of $W(\varphi,m)$ as elements of
$\wreath{S_m}{S_{C_m}}$. Once the generators of $W(\varphi,m)=G$ are
determined, the orbit $O_0$ and the stabilizer $P_0$ of the bracketing
labelled $0$ are returned by
\begin{displaymath}
    \text{\texttt{BlockStabilizerAction(}}G,\ \text{\texttt{[1..m])}}.
\end{displaymath}
The batch function
\begin{displaymath}
    \text{\texttt{IsNice(}}m,\ \varphi\text{\texttt{)}}
\end{displaymath}
first calculates $O_0$, $P_0$ and then returns true if and only if
$O_0=S_{C_m}$, $P_0=S_m$, i.e., if and only if the identity $\varphi$ is
$m$AC-nice.

\begin{example} Here is a transcript of the GAP calculations. The results were
obtained almost instantaneously. (This will not be true for larger values of
$m$.)

\texttt{gap> G := Group( GeneratorsByIdentity( 5, "(xy)z=y(zx)") );}

\texttt{<permutation group with 42 generators>}

\texttt{gap> Size(G);}

\texttt{5596490888974887121059840000000000000000}

\texttt{gap> IsNice( 5, "(xy)z=y(zx)");}

\texttt{true}\\
It is worth noting that $G$ is not all of the wreath product
$\wreath{S_5}{S_{C_5}}$; it is a subgroup of index $2$.
\end{example}

\section{A Short Equational Basis for Boolean Algebras}\label{Sc:BA}

\noindent Several authors have observed that a quasigroup satisfying
$(\ref{Eq:H})$ is an abelian group. (The earliest reference appears to be
\cite{Ho}.) We give a direct and more general proof of this fact based on
ultimate AC-niceness of $(\ref{Eq:H})$.

We say that a groupoid has \emph{one-sided cancellation} if either
(i) $xy=xz$ implies $y=z$ for every $x$, $y$, $z$, or (ii) $yx=zx$
implies $y=z$ for every $x$, $y$, $z$.

\begin{lemma}\label{Lm:AG}
Let $G$ be a groupoid satisfying $\eqref{Eq:H}$. If
$G$ has a neutral element or if $G$ has one-sided cancellation
then it is commutative and associative.
\end{lemma}
\begin{proof}
Consider a product $u$ consisting of $m<5$ elements of $G$. When $G$ has a
neutral element $1$, we can extend $u$ to a product of $5$ elements by
letting $v=((u\cdot 1)\cdot 1)\cdots 1$. When $G$ has one-sided cancellation,
say cancellation on the right, we can pick an element $g\in G$ and extend $u$
into a product of $5$ elements by letting $w=((u\cdot g)\cdot g)\cdots g$.

When $u=xy$, let $u'=yx$. When $u=x(yz)$, let $u'=(xy)z$. (We only discuss
these two cases since we are only interested in commutativity and
associativity.) Let $v'$ (resp.\ $w'$) be the product $v$ (resp.\ $w$) in which
$u$ is replaced by $u'$. Since every groupoid satisfying $(\ref{Eq:H})$ is
$5$AC-nice \cite{HJM1993}, we conclude that $v=v'$, $w=w'$. But $((u\cdot
1)\cdot 1)\cdots 1 = v=v' =((u'\cdot 1)\cdot 1)\cdots 1$ yields $u=u'$ because
$1$ is a neutral element, and, similarly, $w=w'$ yields $u=u'$ because $g$ can
be cancelled on the right.
\end{proof}

We conclude this paper with an application of $(\ref{Eq:H})$ to Boolean
algebras.

Finding short equational bases for varieties of algebras is an important
project in algebra. The variety of Boolean algebras has traditionally occupied
a privileged position in this regard. As early as 1933, Huntington
\cite{Hu1,Hu2} showed that the following three equations form an appealing
short basis for the variety of Boolean algebras:
\begin{align}
    &x + y = y + x,\label{Eq:Hu1}\\
    &(x + y) + z = x + (y + z),\label{Eq:Hu2}\\
    &n(n(x) + y) + n(n(x) + n(y)) = x.\label{Eq:Hu3}
\end{align}
Shortly thereafter, Robbins conjectured that \eqref{Eq:Hu3} could be replaced
with the following shorter equation, which has since come to be known as the
\emph{Robbins equation}:
\begin{equation}\label{Eq:Robbins}
     n(n(x + y) + n(x + n(y))) = x.
\end{equation}
That is, he conjectured that \eqref{Eq:Hu1} and \eqref{Eq:Hu2}, together with
the Robbins equation form an even shorter basis for the variety of Boolean
algebras. But a proof remained elusive for nearly 70 years. The \emph{Robbins
Problem}, as it came to be known, was one of the celebrated open problems in
algebra for most of the 20th century. It was one of Tarski's favorite problems
\cite{Mu}.

Finally, in 1997, Bill McCune solved the problem using his automated theorem
prover, OTTER \cite{Mc}. The buzz generated by McCune's accomplishment was loud
enough to warrant coverage in the New York Times \cite{NYT}! We use equation
\eqref{Eq:H} to offer an even shorter basis for Boolean algebras (Theorem
\ref{Th:BA}).

\begin{theorem}\label{Th:BA}
The following two equations form a basis for the variety of
Boolean algebras:
\begin{align}
    &(x + y) + z = y + (z + x),\tag{additive version of (\ref{Eq:H})}\\
    &n(n(x + y) + n(x + n(y))) = x.\tag{Robbins equation}
\end{align}
\end{theorem}
\begin{proof}
We offer a computer generated proof, found by OTTER \cite{Mc}, that the two
identities imply $x+y=y+x$. Associativity of $+$ then follows. For a primer
on OTTER proofs see \cite{Mc} or \cite{Ph}.

\begin{small}\begin{verbatim}
2 [] (x*y)*z=y* (z*x).
3 [] n(n(x*y)*n(x*n(y)))=x.
5 [] A*B!=B*A.
6 [copy,5,flip.1] B*A!=A*B.
7 [copy,2,flip.1] x* (y*z)= (z*x)*y.
8 [para_into,2.1.1.1,2.1.1] (x* (y*z))*u=y* (u* (z*x)).
9 [copy,8,flip.1] x* (y* (z*u))= (u* (x*z))*y.
10 [para_into,7.1.1.2,7.1.1] x* ((y*z)*u)= ((u*y)*x)*z.
11 [para_into,7.1.1.2,2.1.1] x* (y* (z*u))= (z*x)* (u*y).
12 [para_into,7.1.1,2.1.1] x* ((y*z)*u)= (z* (u*x))*y.
18 [para_into,3.1.1.1.1.1,7.1.1] n(n((x*y)*z)*n(y*n(z*x)))=y.
20 [para_into,3.1.1.1.1.1,2.1.1] n(n(x* (y*z))*n((z*x)*n(y)))=z*x.
22 [para_into,3.1.1.1.1,3.1.1] n(x*n(n(x*y)*n(n(x*n(y)))))=n(x*y).
26 [para_into,3.1.1.1.2.1,2.1.1] n(n((x*y)*z)*n(y* (n(z)*x)))=x*y.
36 [para_into,8.1.1,2.1.1] (x*y)* (z*u)=x* (z* (y*u)).
43 [copy,36,flip.1] x* (y* (z*u))= (x*z)* (y*u).
57 [para_into,9.1.1.2,7.1.1] x* ((y*z)*u)= (y* (x*u))*z.
67 [copy,57,flip.1] (x* (y*z))*u=y* ((x*u)*z).
85 [para_into,10.1.1.2.1,7.1.1] x* (((y*z)*u)*v)= ((v*z)*x)* (u*y).
93,92 [para_into,10.1.1,7.1.1,flip.1] ((x*y)*z)*u= (x*z)* (y*u).
109 [back_demod,85,demod,93,93] x* ((y*z)* (u*v))= (v*x)* (u* (z*y)).
146 [para_from,11.1.1,9.1.1.2,demod,93]
    x* ((y*z)* (u*v))=(y* (x*v))* (u*z).
153 [copy,146,flip.1] (x* (y*z))* (u*v)=y* ((x*v)* (u*z)).
169,168 [para_into,12.1.1.2,2.1.1,flip.1] (x* (y*z))*u=z* (x* (y*u)).
187 [back_demod,153,demod,169] x* (y* (z* (u*v)))=z* ((y*v)* (u*x)).
238 [back_demod,67,demod,169] x* (y* (z*u))=z* ((y*u)*x).
320 [para_into,43.1.1.2,7.1.1] x* ((y*z)*u)= (x*u)* (z*y).
362 [para_into,18.1.1.1.1.1.1,11.1.1,demod,169,169]
    n(n(x* ((y*z)* (u*v)))*n(u* (x*(y*n(v*z)))))=x* (y*u).
397,396 [para_from,18.1.1,3.1.1.1.1]
    n(x*n(n((y*x)*z)*n(n(x*n(z*y)))))=n((y*x)*z).
426 [para_from,92.1.1,18.1.1.1.2.1,demod,169]
    n(n(x* (y* ((z*u)*v)))*n((z*x)*(u*n(v*y))))= (z*u)*x.
529,528 [para_from,168.1.1,3.1.1.1.2.1,demod,169]
    n(n(x*(y* (z*u)))*n(x* (y* (z*n(u)))))=y* (z*x).
562 [para_into,20.1.1.1.1.1,92.1.1,demod,169]
    n(n((x*y)* (z* (u*v)))*n(y* (v*((x*z)*n(u)))))=v* ((x*z)*y).
580 [para_into,20.1.1.1.2.1.1,11.1.1,demod,169,169]
    n(n(x* (y* (z* (u*v))))*n(y*((z*v)* (x*n(u)))))=v* (y* (z*x)).
592 [para_into,20.1.1.1.2.1,2.1.1] n(n(x*(y*z))*n(x* (n(y)*z)))=z*x.
621 [para_into,238.1.1.2,7.1.1] x* ((y*z)*u)=u*((z*y)*x).
663 [para_from,238.1.1,20.1.1.1.1.1.2,demod,169,169,169]
    n(n(x* (y*((z*u)*v)))*n((y*x)* (u* (z*n(v)))))=u* (z* (y*x)).
854 [para_into,22.1.1.1.2.1.1.1,7.1.1,demod,397]
    n((x*y)*z)=n(y* (z*x)).
913 [para_into,854.1.1.1,168.1.1] n(x* (y* (z*u)))=n((z*x)* (u*y)).
1552 [para_into,26.1.1.1.2.1.2,238.1.1,demod,169,169,169]
    n(n((x*y)* (z*(u*v)))*n(y* (x* ((u*z)*n(v)))))=z* (u* (x*y)).
1677 [para_from,621.1.1,7.1.1.2,demod,93,169]
    x* (y* ((z*u)*v))=x* ((u*y)* (z*v)).
1713 [copy,1677,flip.1] x* ((y*z)* (u*v))=x* (z* ((u*y)*v)).
4027,4026 [para_into,109.1.1,320.1.1,demod,169,flip.1]
    (x*y)* (z* (u*v))=x* (y* (z*(u*v))).
5558 [back_demod,1552,demod,4027]
    n(n(x* (y* (z* (u*v))))*n(y* (x*((u*z)*n(v)))))=z* (u* (x*y)).
5895 [back_demod,663,demod,4027]
    n(n(x* (y*((z*u)*v)))*n(y* (x* (u* (z*n(v))))))=u* (z* (y*x)).
5929 [back_demod,562,demod,4027]
    n(n(x* (y* (z* (u*v))))*n(y* (v* ((x*z)*n(u)))))=v*((x*z)*y).
8313,8312 [para_into,187.1.1,238.1.1,demod,169,flip.1]
    x* ((y*z)*(u*v))=x* (z* (y* (u*v))).
10221,10220 [back_demod,1713,demod,8313,flip.1]
    x*(y* ((z*u)*v))=x* (y* (u* (z*v))).
10325,10324 [back_demod,580,demod,8313]
    n(n(x* (y* (z* (u*v))))*n(y* (v* (z* (x*n(u))))))=v* (y* (z*x)).
10345,10344 [back_demod,362,demod,8313]
    n(n(x* (y* (z* (u*v))))*n(u* (x* (z*n(v*y)))))=x*(z*u).
11602 [back_demod,5929,demod,10221,10325]
    x* (y* (z*u))=x* ((u*z)*y).
11604,11603 [back_demod,5895,demod,10221]
    n(n(x* (y* (z* (u*v))))*n(y* (x* (z*(u*n(v))))))=z* (u* (y*x)).
11622 [back_demod,5558,demod,10221,11604]
    x* (y*(z*u))=x* (y* (u*z)).
11820 [back_demod,426,demod,10221]
    n(n(x* (y* (z*(u*v))))*n((u*x)* (z*n(v*y))))= (u*z)*x.
11855 [para_from,11602.1.1,913.1.1.1]
    n(x* ((y*z)*u))=n((z*x)* (y*u)).
11865 [para_from,11602.1.1,20.1.1.1.1.1,demod,93]
    n(n(x* ((y*z)*u))*n((z*x)*(y*n(u))))= (z*y)*x.
11883 [copy,11855,flip.1] n((x*y)* (z*u))=n(y* ((z*x)*u)).
12163,12162 [para_into,11622.1.1.2,7.1.1]
    x* ((y*z)*u)=x* (z* (y*u)).
12176,12175 [back_demod,11883,demod,12163]
    n((x*y)* (z*u))=n(y* (x* (z*u))).
12190,12189 [back_demod,11865,demod,12163,12176,529,flip.1]
    (x*y)*z=x* (y*z).
12505 [back_demod,11820,demod,12190,10345,12190]
    x* (y*z)=z* (y*x).
14989,14988 [para_into,592.1.1.1.1.1,12505.1.1]
    n(n(x* (y*z))*n(z* (n(y)*x)))=x*z.
14996 [para_into,592.1.1.1.2.1,12505.1.1,demod,14989] x*y=y*x.
14997 [binary,14996.1,6.1] $F.
\end{verbatim}\end{small}
\end{proof}

\begin{remark} It is tempting to try to apply Lemma \ref{Lm:AG} in order to
prove Theorem \ref{Th:BA}. Unfortunately, this is only possible if one shows
that there is a neutral element $0$ with respect to addition, necessarily
equal to $n(x+n(x))$ for any $x$. (The addition $+$ is not cancellative.) We
were unable to prove the existence of $0$ without first establishing
commutativity of addition.
\end{remark}

\section{Acknowledgement}

\noindent We thank Alexander Hulpke of Colorado State University for his help
with implementation of wreath products arising from linear identities in GAP.
We also thank Jaroslav Je\v{z}ek of Charles University for bringing several
papers on linear groupoids to our attention.


\begin{small}

\bibliographystyle{plain}

\end{small}

\end{document}